\documentclass[11pt]{article}
\usepackage{amsmath}
\usepackage{amssymb}
\newtheorem{theorem}{Theorem}[section]   

\newtheorem{proposition}[theorem]{Proposition}

\newtheorem{lemma}[theorem]{Lemma} 
\newtheorem{example}[theorem]{Example}

\newtheorem{algorithm}[theorem]{Algorithm}
\def\qed{\hfill $\Box$\medskip}

\def\IC{{\mathbb C}}

\begin{document} 
\openup 0.96 \jot

\title{Products of positive semi-definite matrices}
\author{Jianlian Cui\thanks{
Department of Mathematics, Tsinghua University, 
Beijing 100084, P.R.China.
(jcui@math.tsinghua.edu.cn)} \quad
Chi-Kwong Li\thanks{Department of Mathematics, College of William and Mary, Williamsburg, VA 23187, USA.
(ckli@math.wm.edu)} \quad
Nung-Sing Sze\thanks{Department of Applied Mathematics, the Hong Kong Polytechnic University. (raymond.sze@polyu.edu.hk)}
}
\date{}
\maketitle

\centerline{\bf Dedicated to Professor Rajendra Bhatia.}

\begin{abstract}
It is known that every complex square matrix with nonnegative determinant is the product
of positive semi-definite matrices. There are characterizations of matrices that 
require two or five positive semi-definite matrices in the product.
However, the characterizations of matrices that require three or four positive semi-definite
matrices in the product are lacking. In this paper, we give a complete characterization 
of these two types of matrices. With these results, we give an algorithm 
to determine whether a square matrix can be expressed
as the product of $k$  positive semi-definite matrices but not fewer, for 
$k = 1,2,3,4,5$.
\end{abstract}



\section{Introduction}

Let $M_n$ be the set of $n\times n$ complex matrices.
In \cite{Wu}, the author
showed that a matrix in $M_n$ with nonnegative determinant can always be
written as the product of five or fewer positive semi-definite
matrices. This is an extension to the result in \cite{B} asserting that
every matrix in $M_n$ with positive determinant is the product of 
five or fewer positive definite matrices.  Analogous to the analysis in \cite{B},
the author of \cite{Wu} studied those matrices which can be
expressed as the product of two, three, or four positive semi-definite matrices.
In particular, characterization was obtained for the matrices that can be expressed as
the product of two semi-definite matrices; also, it was shown  that any matrices not of the 
form $zI$, where $z$ is not a nonnegative number, could be written as the product of four
positive semi-definite matrices.
Moreover, it was proved in \cite[Theorem 3.3]{Wu} that if one applies a unitary
similarity to change the square matrix to the form
$T = \left[\begin{matrix}T_1 & * \cr 0 & T_2 \cr\end{matrix}\right] \in M_n$
such that $T_1$ is invertible and  $T_2$ is  nilpotent, and   
if $T_1$ is the product of 
three positive definite matrices, then $T$ can be expressed as the
product of three positive semi-definite matrices as well.
It was suspected that the converse of this statement
is also true. However, the following example shows that it is not the case.

\begin{example}
Let $T = \left[\begin{matrix} -9 & -9 \cr ~0 &~0 \cr\end{matrix}\right]$. Then 
$T =  \left[\begin{matrix} 9 & 3 \cr 3 & 2 \cr\end{matrix}\right]
\left[\begin{matrix} 13 & -15 \cr -15 & 18 \cr\end{matrix}\right]
\left[\begin{matrix} 1 & 1 \cr 1 & 1 \cr\end{matrix}\right]$
is the product of three positive semi-definite matrices.
However, $T_1 = [-9]$ is not the 
product of three positive definite matrices because $\det(T_1) < 0$.
\end{example}

Of course, one may impose the obvious necessary condition that 
$\det(T_1) > 0$, and ask whether the conjecture is valid with this additional 
assumption. Nevertheless, it is easy to modify Example 1 by 
considering
$$T \otimes I_2 =  \left[\begin{matrix} -9I_2 & -9I_2 \cr ~0_2 &~0_2 \cr\end{matrix}\right]
= (A_1\otimes I_2)(A_2 \otimes I_2)(A_3\otimes I_3)$$
using the factorization $T = A_1 A_2 A_3$.
In the modified example, we have $T_1 = -9I_2$ and $T_2 = 0_2$.
By \cite[Theorem 4]{B},  
$T_1 = -9I_2$ is the product of no fewer than five positive definite matrices.
 
In the next section, we will give a complete characterization for those matrices
that can be written as the product of three positive semi-definite matrices and not fewer.
Also, we add an easy to check necessary and sufficient
condition for invertible matrices that can be written as the product of three 
positive definite matrices. With these results, one can use the Jordan form
of a given matrix and its numerical range to decide whether it can be expressed
as the product of two, three, four or five positive semi-definite matrices.

\section{Product of three positive semi-definite matrices}

We will prove the following.

\begin{theorem} \label{main} Suppose
\begin{equation}\label{T-form}
T = \left[\begin{matrix}T_1 & R \cr 0 & T_2 \cr\end{matrix}\right]
\quad \hbox{ such that } T_1  \hbox{ is invertible and } T_2  \hbox{ is nilpotent}.
\end{equation}
Then $T$ is a product of three positive semi-definite matrices if and only if
one of the following holds.
\begin{itemize}
\item[{\rm (a)}] $R \ne 0$ or $T_2 \ne 0$.
\item[{\rm (b)}] $R = 0$, $T_2 = 0$, and $T_1$ is the product 
of three positive definite matrices.
\end{itemize}
\end{theorem}

We establish some lemmas to prove Theorem \ref{main}. The first one is 
covered by \cite[Theorem 1]{B}.

\begin{lemma} \label{SAS} Let $A, S \in M_n$ such that $S$ is invertible. Then 
$A$ is a product of an odd number of  positive semi-definite matrices if and only if
$S^*AS$ is. 
\end{lemma}

In the next lemma, we need the concept of the numerical range of a matrix $A \in M_n$
defined by
$$W(A) = \{ x^*Ax: x\in \IC^n, x^*x = 1\}.$$
The numerical range is a useful tool in studying matrices. 
One may see \cite[Chapter 1]{HJ} for the basic properties of the numerical range.

\begin{lemma} \label{RS}
Suppose $T = \left[\begin{matrix}T_1 & R \cr 0 & 0_p \cr\end{matrix}\right]$
such that $T_1 \in M_m$ is invertible and $R$ is nonzero.
Then $T$ is a product of three positive semi-definite matrices.
\end{lemma}

\it Proof. \rm First, we show that there is a $p\times m$  matrix $S$ 
such that $T_1 + RS$ is the product of three positive definite  matrices.

If $m = 1$, there is  $S \in \IC^p$ such that $T_1 + RS > 0$.
Suppose $m > 1$ and $R$ has a singular value decomposition 
$R = s_1 x_1y_1^* + \cdots + s_k x_k y_k^*$,
where $s_1, \dots, s_k$ are nonzero singular values of $R$,
and  $x_1,\dots,x_k \in \IC^m$ and $y_1,\dots,y_k\in \IC^p$ 
are their corresponding right and left unit singular vectors accordingly.
Let $\{e_1,\dots,e_m\}$ be the standard basis for $\IC^m$.
Take a unitary $U$ such that $U x_1 = e_1$.
Since $T_1$ is invertible, so is $\hat T_1 = UT_1U^*$.
Suppose $\hat t_1^*$ is the first row of $\hat T_1$.
Let $v = \hat t_1 + \epsilon e_2$ with $\epsilon > 0$ and 
$\hat T_1(\epsilon)$ be the $m \times m$ matrix
obtained from $\hat T_1$ by replacing its first row with $v^*$.
Take a sufficiently small $\epsilon  >0$  such that $v$ is not a multiple of $e_1$ and  
the matrix $\hat T_1(\epsilon)$ is still invertible.

Set $S = s_1^{-1} r e^{i\theta} y_1 v^*U$ with $r > 0$ and $\theta\in [0,2\pi)$. Then
$$U(T_1 + RS)U^* = UT_1U^* + URSU^* = \hat T_1 + r e^{i\theta} e_1 v^*.$$
Since $v$ is not a multiple of $e_1$, the rank one matrix $e_1v^*$ is not normal and so 
$W(e_1v^*)$ is an elliptical disk with foci $0$ and $v_1$ 
with length of minor axis $\sqrt{\|v\|^2 - |v_1|^2} > 0$,
where $v_1$ is the first entry of $v$; for example, see \cite[Theorem 1.3.6]{HJ}. 
By the fact that the map $X \mapsto W(X)$ is continuous,
for a sufficiently large $r >0$, $W\left( e_1 v^* + {e^{-i\theta} } \hat T_1 / r \right)$ 
still contains $0$ as an interior point for any $\theta \in [0,2\pi)$. Then so does
$$W\left(\hat T_1 + r e^{i\theta} e_1 v^*\right) 
= re^{i\theta}\, W\left( e_1 v^* + \frac{e^{-i\theta} }{r} \hat T_1 \right).$$ 
In addition, the value $r$ can be chosen so that 
$r > |\det(\hat T_1)/\det(\hat T_1(\epsilon))|$.
Now by the linearity of determinant with respect to the first row, 
$$\det(T_1 + RS) = \det(U(T_1+RS)U^*) = \det(\hat T_1 + 
r e^{i\theta} e_1 v^*) = \det(\hat T_1) + r e^{i\theta} \det (\hat T_1(\epsilon)).$$
Since $|\det(\hat T_1)| < | r e^{i\theta} \det (\hat T_1(\epsilon))|$,  there is 
$\theta\in [0,2\pi)$ such that $\det(T_1 + RS)  > 0$. 
By \cite[Theorem 3]{B} (see also Proposition \ref{three}), $T_1+RS$ is a product
of three positive definite matrices.

Finally, note that 
$$\tilde T = 
\begin{bmatrix}
I_m & S^* \cr
0 & I_p 
\end{bmatrix}
\begin{bmatrix}
T_1 & R \cr 
0 & 0_p 
\end{bmatrix}
\begin{bmatrix}
I_m & 0 \cr
S & I_p 
\end{bmatrix}
= \begin{bmatrix}
T_1 + RS & R \cr 
0 & 0_p 
\end{bmatrix}.
$$
By  \cite[Theorem 3.3]{Wu}, $\tilde T$ is a product of 
three positive semi-definite matrices, and so is $T$ by Lemma \ref{SAS}.
\qed

\medskip
\it Proof of Theorem \ref{main}. \rm
Suppose $T$ is the product of three positive semi-definite matrices.
If $R \ne 0$ or $T_2 \ne 0$, then we are done. Else, $T = T_1 \oplus 0_p$
is the product of three positive semi-definite matrices. By 
\cite[Proposition 3.5]{Wu}, $T_1$ is the product of three positive definite matrices.

\medskip
To prove the converse, we consider the following three cases.

\noindent
{\bf Case 1.} Suppose $R \ne 0$. We use induction on $p$, the size of $T_2$.
If $p = 1$, the result follows from Lemma \ref{RS} as $T_2 = [0]$.
Assume the result holds for $T_2$ with size at most $p-1$.
Since $T_2$ is nilpotent, without loss of generality, we may assume that $T_2$ 
is upper triangular with 
$T_2 = \begin{bmatrix}
T_{21} & T_{22} \cr
0 & 0 \end{bmatrix}$, where $T_{21} \in M_{p-1}$.
Write $R = \begin{bmatrix} R_1 & R_2 \end{bmatrix}$ where $R_1$ is $m \times (p-1)$.
If $R_1 \ne 0$, then by induction, the matrix $\begin{bmatrix} T_1 & R_1 \cr 0 & T_{21} \end{bmatrix}$
is a product of three positive semi-definite matrices, say, $P_1P_2P_3$. Further, 
by  \cite[Theorem 2.2]{Wu} (see also Proposition \ref{two}),
we may assume that both $P_1$  and $P_2$ are invertible. 
Let $X = \begin{bmatrix} R_2 \cr T_{22} \end{bmatrix}$
with size $(m+p-1) \times 1$.
Then for any $\epsilon > 0$,
$$\begin{bmatrix}
P_1 & 0 \cr
0 & 0 
\end{bmatrix}
\begin{bmatrix}
P_2 & \epsilon P_1^{-1} X \cr
\epsilon (P_1^{-1} X)^* & 1 
\end{bmatrix}
\begin{bmatrix}
P_3 & 0 \cr
0 & \epsilon^{-1}
\end{bmatrix}
= \begin{bmatrix}
P_1 P_2 P_3 & X \cr
0 & 0
\end{bmatrix}
= \begin{bmatrix}
T_1 & R_1 & R_2 \cr
0 & T_{21} & T_{22} \cr
0 & 0 & 0 
\end{bmatrix} = T.
$$
Clearly, $Q_1 = P_1 \oplus [0]$
and $Q_3 = P_3 \oplus [\epsilon^{-1}]$ are positive semi-definite matrices. 
Now one can choose  a sufficiently small $\epsilon > 0$ so that
$Q_2 =
\begin{bmatrix}
P_2 & \epsilon P_1^{-1}X \cr
\epsilon (P_1^{-1} X)^* & 1 
\end{bmatrix}$
is also positive semi-definite.
Thus, $T$ is a product of three positive semi-definite matrices $Q_1Q_2Q_3$. 

\medskip
Now suppose $R_1 = 0$. Then the $(m+1)$th column of $T$ is a zero column. By interchanging 
the $(m+1)$th and the last indices, one can see that 
$T$ is permutationally similar to
$$\begin{bmatrix}
T_1 & \tilde R & 0 \cr
0  & \tilde T_{21} & 0 \cr
0 & \tilde T_{22} & 0
\end{bmatrix}
\quad\hbox{where $\tilde R$ is $m \times (p-1)$, $\tilde T_{21}$ is $(p-1) \times (p-1)$, 
and $\tilde T_{22}$ 
is $1 \times (p-1)$.}$$ 
Notice also that 
$\tilde R$ is nonzero and $\tilde T_{21}$ is nilpotent. By induction, 
$\begin{bmatrix}
T_1 & \tilde R  \cr
0  & \tilde T_{21}
\end{bmatrix}$ is a product of three positive semi-definite matrices $P_1P_2P_3$.
By \cite[Theorem 2.2]{Wu},
we can further assume that both $P_2$ and $P_3$ are invertible. 
Let $Y = \begin{bmatrix} 0 & \tilde T_{22} \end{bmatrix}$ with size $1 \times (m+p-1)$. 
Then for any $\epsilon > 0$,
$$\begin{bmatrix}
P_1 & 0 \cr
0 & \epsilon^{-1}
\end{bmatrix}
\begin{bmatrix}
P_2 & \epsilon (Y P_3^{-1})^* \cr
\epsilon (Y P_3^{-1}) & 1 
\end{bmatrix}
\begin{bmatrix}
P_3 & 0 \cr
0 & 0
\end{bmatrix}
= \begin{bmatrix}
P_1 P_2 P_3 & 0 \cr
Y & 0
\end{bmatrix}
= \begin{bmatrix}
T_1 & \tilde R & 0 \cr
0  & \tilde T_{21} & 0 \cr
0 & \tilde T_{22} & 0
\end{bmatrix}.$$
Again, one can choose a sufficiently small $\epsilon > 0$ such that 
all three matrices in the left side of the above equation are positive semi-definite.
Thus, $T$ is permutationally similar to a product of three positive semi-definite matrices,
and hence $T$ can also be written as a product of three positive semi-definite matrices.

\medskip
\noindent
{\bf Case 2.} Suppose $R = 0$ and $T_2$ is nonzero. Without loss of generality, we may assume that
$T_2$ is upper triangular with zero diagonal entries while the first row of $T_2$ is nonzero.
Let $Z$ be the $p \times m$ matrix with 1 at the $(1,1)$-entry and zero elsewise.
Then $Z^*T_2 \ne 0$ and $T_2 Z = 0$. Let $S = \begin{bmatrix}
I_m & 0 \cr
Z & I_p 
\end{bmatrix}$. Then
$$
S^*TS =
\begin{bmatrix}
I_m & Z^* \cr
0 & I_p 
\end{bmatrix}
\begin{bmatrix}
T_1 & 0 \cr
0 & T_2 
\end{bmatrix}
\begin{bmatrix}
I_m & 0 \cr
Z & I_p 
\end{bmatrix}
= \begin{bmatrix}
T_1 + Z^*T_2Z & Z^*T_2 \cr
T_2Z&  T_2
\end{bmatrix}
= \begin{bmatrix}
T_1  & Z^*T_2 \cr
0 &  T_2
\end{bmatrix}.$$
Since $Z^*T_2 \ne 0$, the result follows from Case 1 and Lemma \ref{SAS}.

\medskip
\noindent
{\bf Case 3.} Suppose $R = 0$ and $T_2 = 0$. 
If $T_1$ is the product of three positive definite matrices,
then $T = T_1 \oplus 0_p$ is the product of three positive semi-definite matrices.
\qed

Note that Theorem \ref{main} depends on checking an invertible matrix 
is the product of three positive definite matrices. Such conditions 
are given in \cite[Theorem 3]{B}. We restated the result
in terms of the numerical range in the following proposition,
which is based on the discussion in \cite[Theorem 3 and Fact 3.2]{B}.

\begin{proposition} \label{three}
Let $T \in M_n$ be 
such that $\det(T) > 0$. Then $T$ is the product of three positive definite matrices
if and only if one of the following holds.
\begin{itemize}
\item [{\rm (a)}]  $W(T)$ contains 0 as an interior point.
\item [{\rm (b)}] $W(T)$ contains a positive number, and the 
arguments of the eigenvalues of $T$  can be arranged as:
$-\pi < \theta_1 \le \cdots \le \theta_n < \pi$ such that $\sum_{j=1}^n \theta_j = 0$.
 \end{itemize}
 \end{proposition}

Note that in \cite[p.88]{B}, the author required in condition (3.6b), corresponding to our
condition (b), that all real eigenvalues of $T$ are positive, which is ensured by 
our assumption that $\theta_j \in (-\pi, \pi)$ for all $j$.

\section{Determining the number of factors}

In this section, we describe an algorithm 
to determine the smallest number of positive semi-definite matrices whose product
equals a given $A \in M_n$ with 
nonnegative determinant.

\medskip
We first present the following theorem providing some easy tests for a matrix $A$ to be the 
product of two positive semi-definite matrices.
The equivalence of conditions (a), (b), (c) were given in \cite[Theorem 2.2]{Wu}. We include a 
short proof, which is different from that of Wu.

\begin{proposition} \label{two}
Let $A$ be a square matrix. The following are equivalent.
\begin{itemize}
\item[{\rm (a)}] $A$ is the product of two positive semi-definite matrices.
\item[{\rm (b)}] 
$A = BC$, where $B, C$ are positive semi-definite matrices such that 
$B$ or $C$ can be assumed to be invertible.
\item[{\rm (c)}] $A$ is similar to a nonnegative diagonal matrix.
\item[{\rm (d)}] $A$ is unitarily similar to an upper block triangular matrix
such that the diagonal blocks are scalar matrices corresponding to distinct scalars. 
\item[{\rm (e)}] The minimal polynomial of $A$ only has simple nonnegative zeros.
\end{itemize}
\end{proposition}

\it Proof. \rm The equivalence of (c), (d), (e) are clear.

(c) $\Rightarrow$ (b): 
$A = S^{-1}DS 
= S^{-1}(S^{-1})^*(S^*DS) 
= S^{-1}D (S^{-1})^*(S^*S)$.

\medskip
(b) $\Rightarrow$ (a): Trivial.

\medskip
(a) $\Rightarrow$ (c):  Suppose $A = BC$, where $B$ and $C$ are $n\times n$ positive semi-definite.
Let $U$ be unitary such that 
$U^*BU = B_0 \oplus 0_k$, where $B_0 \in M_{n-k}$ is positive definite.
Let  $U^*CU = \left[\begin{matrix}C_{11}& C_{12} \cr C_{21} & C_{22} \cr\end{matrix}\right]$
be such  that $C_{22} \in M_k$.
Assume $V$ is unitary such that $V^*C_{11}V = C_0 \oplus 0_\ell$ 
for a positive definite matrix $C_0\in M_{n-k-\ell}$. 
We may replace $U$ by $U(V\oplus I)$ and assume that
$C_{11} = C_0 \oplus 0_\ell$.
Because $C$ is positive semi-definite, 
$U^*CU = \left[\begin{matrix}C_0 & 0 & C_1  \cr 
0 & 0_\ell & 0  \cr C_1^* & 0 & C_{22} \cr \end{matrix}\right]$.
Then
$$U^*AU = 
\left[\begin{matrix}B_0& 0 \cr 0 & 0_k \cr\end{matrix}\right]
\left[\begin{matrix} C_{0}& 0 &  C_1 \cr 0 & 0_\ell & 0 \cr
C_1^* & 0 & C_{22} \cr\end{matrix}\right]
= \left[\begin{matrix}B_0& 0 \cr 0 & 0_k \cr\end{matrix}\right]
\left[\begin{matrix} C_{0}& 0 &  C_1 \cr 0 & 0_\ell & 0 \cr
0 & 0 & 0_k \cr\end{matrix}\right].$$
Now
\begin{eqnarray*}
&& 
\begin{bmatrix}
B_0^{-\frac{1}{2}} & 0 \cr 
0 & I_k 
\end{bmatrix}
\begin{bmatrix}
I & 0 & C_0^{-1} C_1 \cr 
0 & I & 0  \cr
0 & 0 & I_k
\end{bmatrix}
\left[\begin{matrix}B_0& 0 \cr 0 & 0_k \cr\end{matrix}\right]
\left[\begin{matrix} C_{0}& 0 &  C_1 \cr 0 & 0_\ell & 0 \cr
0 & 0 & 0_k \cr\end{matrix}\right]
\begin{bmatrix}
I & 0 & -C_0^{-1} C_1 \cr 
0 & I & 0  \cr
0 & 0 & I_k
\end{bmatrix}
\begin{bmatrix}
B_0^{\frac{1}{2}} & 0 \cr 
0 & I_k 
\end{bmatrix} \\[1mm]
&=& 
\begin{bmatrix}
B_0^{-\frac{1}{2}} & 0 \cr 
0 & I_k 
\end{bmatrix}
\left[\begin{matrix}B_0& 0 \cr 0 & 0_k \cr\end{matrix}\right]
\left[\begin{matrix} C_{0}& 0 &  0 \cr 0 & 0_\ell & 0 \cr
0 & 0 & 0_k \cr\end{matrix}\right]
\begin{bmatrix}
B_0^{\frac{1}{2}} & 0 \cr 
0 & I_k 
\end{bmatrix}
= \left( B_0^\frac{1}{2} \begin{bmatrix} C_0 & 0 \cr 0 & 0_\ell \end{bmatrix} B_0^{\frac{1}{2}} \right)
\oplus 0_k.
\end{eqnarray*}
Therefore, $A$ is similar to $\left( B_0^\frac{1}{2} \begin{bmatrix} C_0 & 0 \cr 0 & 0_\ell \end{bmatrix} B_0^{\frac{1}{2}} \right)
\oplus 0_k$,
which is positive semi-definite and is (unitarily) similar to a nonnegative diagonal matrix.
\qed

Now, we are ready to present an algorithm to check whether a matrix $A \in M_n$
with $\det(A) \ge 0$ can be written as a product of $k$ positive semi-definite matrices,
but not fewer, for $k = 1, 2, 3, 4, 5$.

\begin{algorithm} \rm Let $A \in M_n$ with $\det(A) \ge 0$.

If $A = \alpha I_n$ such that $\alpha \notin [0, \infty)$, then
$A$ can be expressed as a product of five positive semi-definite matrices, but not fewer.
Otherwise, apply a unitary similarity to $A$ to get an upper triangular matrix 
$$T = \left[\begin{matrix}T_1 & R \cr 0 & T_2 \cr\end{matrix}\right]
\quad \hbox{ such that } T_1  \hbox{ is invertible and } T_2  \hbox{ is  nilpotent}.$$

\begin{itemize}
\item[{\bf (1)}] If $T$ is a nonnegative diagonal matrix, then $A$ is itself a positive
semi-definite matrix.

\item[{\bf (2)}] Condition {\bf (1)} does not hold, and  
$A$ satisfies any one of the equivalent conditions in Proposition \ref{two}.
Then $A$ can be expressed as a product of two positive semi-definite matrices, but not fewer.

\item[{\bf (3)}] Suppose {\bf (1)} and {\bf (2)} do not hold.
Then 
$A$ can be expressed as a product of three positives semi-definite matrices, but not fewer,
if any of the following holds. 
\begin{itemize}
\item[\rm (3.a)] $R$ or $T_2$ is nonzero. 

\item[\rm (3.b)] Both $R = 0$ and $T_2 = 0$, and $T_1$ is the product of 
three positive definite matrices. 
\end{itemize}

\medskip\noindent
In (3.b), the invertible matrix $T_1$ is a product of three positive
definite matrices if one of the following holds.

\begin{itemize}
\item[\rm (i)] $\det(T_1) > 0$ and $W(T_1)$ contains 0 as an interior point,

\item[\rm (ii)] $0$ is not in the interior of $W(T_1)$
and $W(T_1)$ contains a positive number and 
$\sum \theta_j = 0$, where $-\pi < \theta_1 \le  \dots \le \theta_k < \pi$ are 
the arguments of the eigenvalues of $T_1$.
\end{itemize}

\item[{\bf (4)}] 
Suppose conditions {\bf (1), (2), (3)} do not hold, i.e., 
$T = T_1\oplus\, 0_p$ such that neither (i) nor (ii) holds 
for the upper triangular matrix $T_1$.
Then  $A$ can be expressed as a product of four positive semi-definite matrices, but not fewer. 
\end{itemize}
\end{algorithm}

\medskip\noindent
{\bf Acknowledgment}

The authors would like to thank the referee for his/her careful reading of the 
manuscript.
The research of Cui was supported by National Natural Science Foundation of China  11271217.
Li is an affiliated member of the Institute for Quantum Computing, University of 
Waterloo. He is also an honorary professor of the Shanghai University, and 
the University of Hong Kong;  his research was supported by the
USA NSF grant DMS 1331021, the Simons Foundation Grant 351047, and the NNSF of China Grant 11571220.  
The research of Sze was supported by a HK RGC grant PolyU 502512
and a PolyU central research grant G-UC25.

\end{document}